\documentclass[12pt,a4paper]{amsart}
\usepackage[all]{xy}
\usepackage{amscd}
\usepackage{amsfonts}
\usepackage{amsmath}
\usepackage{amssymb}
\usepackage{amsthm}
\usepackage[colorlinks,linkcolor=blue]{hyperref}
\usepackage{calc}
\usepackage{extarrows}
\usepackage{graphicx}
\usepackage{graphics}
\usepackage{mathrsfs}
\usepackage[margin=2.9cm]{geometry}

\theoremstyle{plain}
\newtheorem{thm}{Theorem}[section]
\newtheorem{lema}[thm]{Lemma}
\newtheorem{ppst}[thm]{Proposition}

\theoremstyle{definition}
\newtheorem{defn}[thm]{Definition}
\newtheorem{rem}[thm]{Remark}
\newtheorem{exam}[thm]{Example}

\newcommand{\bigboxplus}{
  \mathop{
    \vphantom{\bigoplus}
    \mathchoice
      {\vcenter{\hbox{\resizebox{\widthof{$\displaystyle\bigoplus$}}{!}{$\boxplus$}}}}
      {\vcenter{\hbox{\resizebox{\widthof{$\bigoplus$}}{!}{$\boxplus$}}}}
      {\vcenter{\hbox{\resizebox{\widthof{$\scriptstyle\oplus$}}{!}{$\boxplus$}}}}
      {\vcenter{\hbox{\resizebox{\widthof{$\scriptscriptstyle\oplus$}}{!}{$\boxplus$}}}}
  }\displaylimits
}

\title[]{A few results on associativity of hypermultiplications in polynomial hyperstructures over hyperfields}
\author{Ziqi Liu}
\address{Jilin University, Changchun, Jilin, China}
\email{liuzq0616@mails.jlu.edu.cn}

\thanks{I thank Matthew Baker for inspiring discussions and valuable suggestions during my visiting to Georgia Institute of Technology. I thank Oliver Lorscheid for helpful comments on my proofs. This work was derived from an unchecked statement in a version of Baker and Lorscheid's paper and the author started this work in October. The author may update this work non-scheduledly. The author was supported in part by Chinese Scholarship Council during his four-month stay in the United States in 2019.
}

\begin{document}

\begin{abstract}
  In Baker and Lorscheid's paper, they introduce a new hyperstructure: the polynomial hyperstructure Poly$(\mathbb{F})$ over a hyperfield $\mathbb{F}$. In this work, the author focuses on associativity of hypermultiplications in those hyperstructures and gives elementary propositions. The author also shows examples of polynomial hyperstructures over hyperfields with non-associative hypermultiplications. Then, he proves that though the hypermultiplication in Poly$(\mathbb{T})$ is associative for linear polynomials, it is not associative in general. Moreover, he shows that if $1\boxplus_{\mathbb{F}}1$ is not a singleton for hyperfield $\mathbb{F}:=(\mathbb{F},\odot,\boxplus_{\mathbb{F}},1,0)$, the hypermultiplication in Poly$(\mathbb{F})$ is not associative.
\end{abstract}

\maketitle
\tableofcontents

\newpage
\section{Hypefields and Polynomials over a Hyperfield}
\subsection{The Definition of Hyperfields}
\begin{defn}
A \textbf{hyperoperation} on a set $S$ is a map $\square:S\times S\rightarrow 2^S\backslash\{\varnothing\}$. Moreover, for a given hyperoperation $\square$ on $S$ and non-empty subsets $A,B$ of $S$, $A\,\square\,B$ is define as
$$A\,\square\,B:=\bigcup_{a\in A,b\in B}(a\,\square\,b)$$
A hyperoperation $\square$ in $S$ is \textbf{commutative} if $a\,\square\,b=b\,\square\,a$ for all $a,b\in S$. If not especially mentioned, hyperoperations in this work will always be commutative.\\
A hyperoperation $\square$ in $S$ is \textbf{associative} if $a\,\square\,(b\,\square\,c)=(a\,\square\,b)\square\,c$ for all $a,b,c\in S$.
\end{defn}
\begin{defn}
Given an associative hyperoperation $\square$ in $S$, a \textbf{hypersum} is recursively defined as
$$x_1\,\square\,\cdots\,\square\,x_n:=\bigcup_{x'\in x_2\,\square\,\cdots\,\square\,x_n}x_1\,\square\,x'$$
for $x_1,\dots,x_n$ in $S$ where $n\geq2$.
\end{defn}
\begin{defn}
A \textbf{hypergroup} is a tuple $(G,\boxplus,0)$, where $\boxplus$ is an associative hyperoperation on $G$ such that:\\
(1) $0\boxplus x=\{x\}$ for all $x\in G$;\\
(2) For every $x\in G$ there is a unique element $-x$ of $G$ such that $0\in x\boxplus-x$;\\
(3) $x\in y\boxplus z$ if and only if $z\in x\boxplus(-y)$.\\
Here $-x$ is often called as the \textbf{hyperinverse} of $x$ and (3) as the reversibility axiom.
\end{defn}
\begin{defn}
A (Krasner) \textbf{hyperring} is a tuple $(R,\odot,\boxplus,1,0)$ such that:\\
(1) $(R,\odot,1)$ is a commutative monoid;\\
(2) $(R,\boxplus,0)$ is a commutative hypergroup;\\
(3) $0\odot x=x\odot0=0$ for all $x\in R$;\\
(4) $a\odot(x\boxplus y)=(a\odot x)\boxplus(a\odot y)$ for all $a,x,y\in R$;\\
(5) $(x\boxplus y)\odot a=(x\odot a)\boxplus(y\odot a)$ for all $a,x,y\in R$.\\
In the following part, we will use the underlying set $R$ to refer to a hyperring and may omit $\odot$ if there is no likehood of confusion. In addition, $1$ is called the \textbf{unit element} and $0$ is called the \textbf{zero element} in the hyperring $R$.
\end{defn}
\begin{defn}
A hyperring $F$ is called a \textbf{hyperfield} if $0\neq 1$ and every non-zero element of $F$ has a multiplicative inverse.
\end{defn}

\begin{exam}
If $(\mathbb{F},\cdot,+)$ is a field, then $\mathbb{F}$ can be trivially associated with a hyperfield $(\mathbb{F},\odot,\boxplus)$ where $x\odot y=x\cdot y$ and $x\boxplus y=\{x+y\}$ for all $x,y\in\mathbb{F}$.\\
In the following context, when we mention a field $\mathbb{F}$, we may actually refer to the hyperfield associated with $\mathbb{F}$.
\end{exam}
\begin{exam}
Consider $\mathbb{K}=(\{0,1\},\odot,\boxplus,1,0)$ with the usual multiplication rule and a hyperaddition $\boxplus$ defined by
$$0\boxplus0=\{0\},\qquad1\boxplus0=0\boxplus1=\{1\},\qquad 1\boxplus1=\{0,1\}$$
then $\mathbb{K}$ is a hyperfield, called the \textbf{Krasner hyperfield}.
\end{exam}
\begin{exam}
Consider $\mathbb{S}=(\{0,1,-1\},\odot,\boxplus,1,0)$ with the usual multiplication rule and a hyperaddition $\boxplus$ generated by
$$x\boxplus x=\{x\},\quad x\boxplus 0=\{x\},\quad1\boxplus-1=\{-1,0,1\}$$
then $\mathbb{S}$ is a hyperfield, called the \textbf{hyperfield of signs}.
\end{exam}

\begin{exam}
Consider $\mathbb{W}=(\{0,1,-1\},\odot,\boxplus,1,0)$ with the usual multiplication rule and a hyperaddition $\boxplus$ generated by
$$x\boxplus x=\{x,-x\},\quad x\boxplus 0=\{x\},\quad1\boxplus-1=\{-1,0,1\}$$
then $\mathbb{W}$ is a hyperfield, called the \textbf{weak hyperfield of signs}.
\end{exam}
\begin{rem}
More generally, given a multiplicatively written abelian group $(G,\cdot,1)$ and a self-inverse element $e$ of $G$, there exists a hyperfield $W(G,e)=(G\cup\{0\},\cdot,\boxplus,1,0)$ where the multiplication $\cdot$ is the same as that in $G$ with $0\cdot x=0$ for all $x\in G\cup\{0\}$, and the hyperaddition is defined by
$$0\boxplus x=\{x\},\quad x\boxplus(e\cdot x)=G\cup\{0\},\quad x\boxplus y=G$$
for any nonzero $x$ and $y$ with $y\neq ex$. Such a hyperfield is called a \textbf{weak hyperfield}.
\end{rem}
\begin{exam}
Let $\mathbb{T}:=\mathbb{R}\cup\{-\infty\}$ as sets and define hyperoperation $\boxplus$ as
$$ x\boxplus y=\left\{
\begin{aligned}
\{\max\{x,y\}\} &, & x\neq y \\
\{z\in\mathbb{T}:z\leq x\} &,  & x=y
\end{aligned}
\right.
$$
and $\odot$ as $x\odot y=x+y$. Then $\mathbb{T}$ is a hyperfield, called the \textbf{tropical hyperfield}.
\end{exam}
\begin{rem} More generally, let $\Gamma$ be a totally ordered abelian group (written multiplicatively) one can define a canonical hyperfield structure on set $\Gamma\cup\{0\}$ where
\begin{itemize}
  \item the multiplication $\odot$ is multiplication in $\Gamma$ with $0\odot x=0$ for all $x\in\Gamma\cup\{0\}$
  \item the hyperaddition $\boxplus$ is defined as $x\boxplus x:=\{y:y\leq x\}$ and $x\boxplus y:=\max\{x,y\}$ for $x\neq y$.
\end{itemize}
and $x\geq0$ for all $x\in\Gamma$. Such a hyperfield is called a \textbf{valuative hyperfield}. In the tropical hyperfield $\mathbb{T}$, the zero element is $-\infty$ and the unit element is $0$. In addition, the Krasner hyperfield $\mathbb{K}$ is also a valuative hyperfield.
\end{rem}
\begin{exam}
Let $\mathbb{P}=S^1\cup\{0\}$, where $S^1=\{z\in\mathbb{C}:|z|=1\}$ is the complex unit circle. Then one can define a hyperfield structure on $\mathbb{P}$ where the multiplication is the usual one in the complex field $\mathbb{C}$ and the hyperaddition is defined as
$$ x\boxplus y=\left\{
\begin{array}{cll}
\{x\}, & & y=0 \\
\{0,x,-x\}, & & y=-x \\
\{\frac{\alpha x+\beta y}{||\alpha x+\beta y||}:\alpha,\beta\in\mathbb{R}_+\}, &  & \textup{otherwise}
\end{array} \right. $$
This hyperfield structure on $S^1\cup\{0\}$ is called the \textbf{phase hyperfield}.
\end{exam}
\begin{exam}
Let $\mathbb{V}$ be the set $\mathbb{R}_{\geq0}=\mathbb{R}_+\cup\{0\}$ of nonnegative real numbers with the usual multiplication in the field $\mathbb{R}$ and the hyperaddition is defined as
$$x\boxplus y=\{z\in\mathbb{R}_{\geq0}:|x-y|\leq z\leq x+y\}$$
Then $\mathbb{V}$ is a hyperfield called the \textbf{Viro hyperfield} (or the triangle hyperfield).
\end{exam}

For more information about the construction of hyperfields, \cite{4} and \cite{5} will be good references.
\subsection{Polynomials over Hyperfields}
\begin{defn}
Given a hyperfield $\mathbb{F}$, a \textbf{polynomial} over $\mathbb{F}$ (or with coefficients in $\mathbb{F}$) is a map $p:\mathbb{F}\rightarrow2^{\mathbb{F}}$ that
\begin{center}
$a\longmapsto c_na^n\boxplus c_{n-1}a^{n-1}\boxplus\cdots\boxplus c_1a\boxplus c_0$
\end{center}
where $\{c_i\}_{i=0}^n\subset\mathbb{F}$ and $c_n$ is not the additive unit (zero element) in $\mathbb{F}$.\\
For such $p$, we denote it by $p(T)=c_nT^n+c_{n-1}T^{n-1}+\cdots+c_1T+c_0$. In addition, the \textbf{degree} of $p$ is defined to be the largest nonnegative integer $n$ such that the coefficient of $T^n$ is nonzero.
\end{defn}
\begin{exam}
For any field $\mathbb{F}$, elements in the polynomial ring $\mathbb{F}[T]$ are polynomials over hyperfield $\mathbb{F}$.
\end{exam}
\begin{exam}
The polynomial $p(T)=1T^3+(-2)$ over the tropical hyperfield $\mathbb{T}$ is exactly represented by
$$p(a)=1a^3\boxplus(-2)=\left\{\begin{aligned}
-2\quad&,\,\,a<-1\\
[-\infty,-2]&,\,\,a=-1\\
1a^3\quad&,\,\,a>-1
\end{aligned}\right.$$
where the order $<$ is the same as natural order in $\mathbb{R}$ and $1a^3$ means $1+3a$ in $\mathbb{R}$.
\end{exam}
\begin{ppst}\label{1.16}
Let $\mathbb{F}$ be a hyperfield, the set of all polynomials over $\mathbb{F}$ is naturally endowed with two hyperoperations\\
(1) $p\boxdot q=\{e_{mn}T^{mn}+\cdots+e_1T+e_0 : e_i=\mathop{\bigboxplus_{\mathbb{F}}}\limits_{k+l=i}c_kd_l\}$;\\
(2) $p\boxplus q=\{e_{k}T^{k}+\cdots+e_2T^2+e_1T+e_0 : e_i=c_i\boxplus_{\mathbb{F}} d_i\}$.\\
for $p(T)=c_nT^n+\cdots+c_1T+c_0$ and $q(T)=d_mT^m+\cdots+d_1T+d_0$ with $k=\max\{m,n\}$.\\
Clearly, $\boxdot$ and $\boxplus$ are commutative since $\boxplus_{\mathbb{F}}$ is commutative.
\end{ppst}
\begin{rem} In this paper, we will call this hyperstructure as \textbf{the polynomial hyperstructure over a hyperfield} $\mathbb{F}$ and denote it by $\textup{Poly}(\mathbb{F})$. In \cite{2}, this hyperstructure is called a polynomial hyperring while it is in fact not a hyperring. In other materials like \cite{1}, it is called a superring or a hyperring of polynomials.
\end{rem}

\begin{defn}
Let $p(T)=c_n T^n+ c_{n-1}T^{n-1}+\cdots+ c_1T+ c_0$ be a polynomial over a hyperfield $\mathbb{F}$, an element $a\in\mathbb{F}$ is called a \textbf{root} of $p$ if and only if either the following equivalent conditions is satisfied:\\
(1) $0\in p(a)=c_n a^n\boxplus c_{n-1}a^{n-1}\boxplus \cdots\boxplus c_1a\boxplus c_0$;\\
(2) there exists elements $d_0,d_1,\dots,d_{n-1}\in\mathbb{F}$ such that
\begin{center}
$c_0=-ad_0,\,\,c_i\in -ad_i\boxplus d_{i-1}$ for $i=1,\dots,n-1$ and $c_n=d_{n-1}$
\end{center}
Notice that here (2) means that $p\in(T-a)\boxdot q$ in Poly$(\mathbb{F})$.
\end{defn}

\begin{defn}
Let $\displaystyle p(T)=c_n T^n+ c_{n-1}T^{n-1}+\cdots+ c_1T+ c_0$ be a polynomial over a hyperfield $\mathbb{F}$, if $a$ is not a root of $p$, set mult$_{a}(p)=0$. If $a$ is a root of $p$, define
\begin{center}
mult$_{a}(p)$ = $1+\max\{\textup{mult}_{a}(q):p\in (T-a)q\}$
\end{center}
as the \textbf{multiplicity} of the root $a$ of $p$. Moreover, for a nonempty set $S$, define
\begin{center}
mult$_{S}(p)$ = $1+\max\{\textup{mult}_{S}(q):p\in (T-a)q$ for some $a\in S\}$
\end{center}
It is clearly that mult$_S(p)\leq\textup{deg}(p)$ for any $S\subset\mathbb{F}$ and polynomial $p$ over $\mathbb{F}$.
\end{defn}

\begin{exam}
Given a polynomial $p(T)=T^3-T$ over the hyperfield of signs $\mathbb{S}$, it is clear that $q(T)=T^2-1$ is the only polynomial in $\textup{Poly}(\mathbb{S})$ such that $p\in T\boxdot q$. Then one can see mult$_0(p)=1$ since $0\in p(0)=\{0\}$ and $0$ is not a root of $q(T)$.
\end{exam}

\begin{exam}
Given a subset $S=[1,+\infty)$ of the Viro hyperfield $\mathbb{V}$ and a polynomial $p(T)=T^2+3T+1$ over $\mathbb{V}$, one can see mult$_S(p)\geq1$ since the inequalities
\begin{center}
$|3a-1|\leq a^2\leq 3a+1$
\end{center}
has solutions in $S$. Notice that $(T+a)\boxdot(T+a)=\{T^2+bT+a^2:b\in a\boxplus a\}$, it is clear that $p\notin (T+a)^2$ for any $a\geq1$. Therefore, mult$_S(p)\leq1$ and then mult$_S(p)=1$. Here one should notice that $a=-a$ for any $a$ in $\mathbb{V}$.
\end{exam}

\section{Associativity of Hypermultiplications in Poly$(\mathbb{F})$}
\subsection{A first glance on associativity of hypermultiplications}
Seen from the definition of the polynomial hyperstructure over a hyperfield, the hyperaddtion in any polynomial hyperstructure over a hyperfield is associative, but associativity of hypermultiplication has not yet been determined. In this part, one can see some first thoughts about associativity of hypermultiplication in Poly$(\mathbb{F})$.

\begin{lema}
Given a hyperfield $\mathbb{F}$, then for each polynomial $p$ over $\mathbb{F}$, there exists a monic polynomial $p_0$ over $\mathbb{F}$ and $a\in\mathbb{F}$ such that $p=a\boxdot p_0$.

\begin{proof}
By definition, each element $p$ in Poly$(\mathbb{F})$ can be represented as
$$p(T)=a_nT^n+a_{n-1}T^{n-1}+\cdots+a_1T+a_0$$
with $a_n\neq0$. It is easy to check that
$$p_0(T)=T^n+a_{n-1}a^{-1}_nT^{n-1}+\cdots+a_1a^{-1}_nT+a_0a^{-1}_n$$
satisfies $p(T)=a_n\boxdot p_0(T)$.
\end{proof}
\end{lema}
\begin{rem}
In addition, the following two propositions Proposition \ref{2.4} and Proposition \ref{2.5} imply that one can only check the associativity of the hypermultiplication for monic polynomials.
\end{rem}

\begin{ppst}\label{2.4}
Given a hyperfield $\mathbb{F}$, one has
$$(a\boxdot p_0)\boxdot(b\boxdot q_0)=(ab)\boxdot(p_0\boxdot q_0)$$
for any $a,b\in\mathbb{F}$ and $p_0,q_0\in\textup{Poly}(\mathbb{F})$.
\end{ppst}

\begin{ppst}\label{2.5}
Given a hyperfield $\mathbb{F}$, if the hypermultiplication in the polynomial hyperstructure over $\mathbb{F}$ is associative for all monic polynomials, then the hypermultiplication in $\textup{Poly}(\mathbb{F})$ is always associative.

\begin{proof}
Given three polynomials $p(T),q(T),r(T)$ in $\textup{Poly}(\mathbb{F})$, there exist $a,b,c$ with
$$p(T)=a\boxdot p_0(T),\,\,q(T)=b\boxdot q_0(T),\,\,r(T)=c\boxdot r_0(T)$$
where $p_0(T),q_0(T),r_0(T)$ are monic. Then, one can see
\begin{align*}
p\boxdot(q\boxdot r)&=(a\boxdot p_0)\boxdot((b\boxdot q_0)\boxdot(c\boxdot r_0))\\
&=(a\boxdot p_0)\boxdot(bc\boxdot (q_0\boxdot r_0))\\
&=abc\boxdot( p_0\boxdot(q_0\boxdot r_0))\\
&=abc\boxdot( q_0\boxdot(p_0\boxdot r_0))\\
&=(b\boxdot q_0)\boxdot(ac\boxdot(p_0\boxdot r_0))\\
&=(b\boxdot q_0)\boxdot((a\boxdot p_0)\boxdot(c\boxdot r_0))\\
&=q\boxdot(p\boxdot r)
\end{align*}
and we are done.
\end{proof}
\end{ppst}

\begin{ppst}
Given a hyperfield $\mathbb{F}$ and two polynomials $p,q$ in $\textup{Poly}(\mathbb{F})$, we have
$$T^n\boxdot(p(T)\boxdot q(T))=p(T)\boxdot(T^n\boxdot q(T))$$
for any $n\in\mathbb{N}$.

\begin{proof}
It is easy to check that
$$T^n\boxdot (a_mT^m+a_{m-1}T^{m-1}+\cdots+a_0)=\{a_{m+n}T^m+a_{m-1}T^{m+n-1}+\cdots+a_0T^n\}$$
for any $r(T)=a_mT^m+a_{m-1}T^{m-1}+\cdots+a_0$ in Poly$(\mathbb{F})$ and $n\in\mathbb{N}$. Then we can see
\begin{align*}
T^n\boxdot(p(T)\boxdot q(T))&=\{T^n\boxdot r(T):r(T)\in p(T)\boxdot q(T)\}\\
&=\{T^n\boxdot r(T):r(T)\in p(T)\boxdot q(T)\}\\
&=\{T^n\boxdot r(T):T^n\boxdot r(T)\subset p(T)\boxdot(T^n\boxdot q(T))\}\\
&=\{r(T):r(T)\in p(T)\boxdot(T^n\boxdot q(T))\}\\
&=p(T)\boxdot(T^n\boxdot q(T))
\end{align*}
and we are done.
\end{proof}
\end{ppst}

\subsection{Some non-associative hypermultiplications}
In this part, we will give examples to show non-associativity of some hypermultiplications in polynomial hyperstructures over hyperfields and try to find some clues behind them.
\begin{ppst}
In the Viro hyperfield $\mathbb{V}$, $0\in x\boxplus y$ if and only if $x=y$.

\begin{proof}
If $0\in x\boxplus y=\{z:|x-y|\leq z\leq x+y\}$, then $|x-y|=0$ and then $x=y$.\\
If $x=y$, then $0\in x\boxplus y=x\boxplus x=\{z:0\leq z\leq 2x\}$.
\end{proof}
\end{ppst}
\begin{rem}
This proposition implies that each $x$ in $\mathbb{V}$ is the hyperinverse of itself.
\end{rem}
\begin{exam}\label{V}
Consider the polynomial $p(T)=T^3+2T^2+11T+6$ over the Viro hyperfield $\mathbb{V}$. One is able to check that
$$\begin{aligned}
p(T)\in(T+2)\boxdot(T^2+4T+3)&\subseteq\{T^3+2\boxplus4T^2+3\boxplus 8T+6\}\\
&=(T+2)\boxdot(T^2+4T+3)\\
&\subseteq(T+2)\boxdot((T+1)\boxdot(T+3))
\end{aligned}$$
Then we claim that $p(T)\notin(T+1)\boxdot((T+2)\boxdot(T+3))$. In fact, we have
$$\begin{aligned}
(T+1)\boxdot((T+2)\boxdot(T+3))
&=\{(T+1)\boxdot(T^2+d_1T+6)\,|\,d_1\in[1,5]\}\\
&=\{T^3+(d_1\boxplus1)T^2+(d_1\boxplus6)T+6\,|\,d_1\in[1,5]\}
\end{aligned}$$
Here one can see that if $11\in d_1\boxplus 6$, $d_1$ must be $5$ and then $2\notin d_1\boxplus 1=[4,6]$. Therefore, the hypermultiplication in Poly$(\mathbb{V})$ is not associative.
\end{exam}
\begin{ppst}
For any $x$ in the phase hyperfield $\mathbb{P}$, one has $x\boxplus x=\{x\}$.
\end{ppst}
\begin{exam}\label{P}
Consider the polynomial $p(T)=T^3-e^{i\frac{\pi}{8}}T^2+e^{i\frac{5\pi}{24}}T-e^{i\frac{\pi}{3}}$ over the phase hyperfield $\mathbb{P}$. It is not difficult to check that $p\in(T-e^{i\frac{\pi}{6}})q$ and mult$_a(q)=2$ where $q(T)=T^2-e^{i\frac{\pi}{12}}T+e^{i\frac{\pi}{6}}$ and $a=e^{i\frac{\pi}{12}}$. But, one can obtain that
$$0\notin e^{i\frac{\pi}{4}}\boxplus(-e^{i\frac{7\pi}{24}})\boxplus e^{i\frac{7\pi}{24}}\boxplus(-e^{i\frac{\pi}{3}})=p(e^{i\frac{\pi}{12}})$$
which means that $e^{i\frac{\pi}{12}}$ is not a root of $p$. Therefore, there does not exist a $r\in\textup{Poly}(\mathbb{P})$ such that $p\in(T-e^{i\frac{\pi}{12}})r$, let alone $p\in(T-e^{i\frac{\pi}{12}})\boxdot((T-e^{i\frac{\pi}{6}})\boxdot(T-e^{i\frac{\pi}{12}}))$.\\
Therefore, the hypermultiplication in Poly$(\mathbb{P})$ is not associative.
\end{exam}

\begin{exam}\label{W}
Consider the polynomial $p(T)=T^3-1$ over the weak hyperfield of signs $\mathbb{W}$, one can see
$$p(T)\in (T-1)\boxdot(T^2+T+1)\subseteq(T-1)\boxdot((T+1)\boxdot(T+1))$$
but $0\notin(-1)\boxplus(-1)=p(-1)$, which implies that
$$p(T)\notin (T+1)\boxdot((T-1)\boxdot(T+1))$$
Therefore, the hypermultiplication in Poly$(\mathbb{W})$ is not associative.
\end{exam}

\begin{defn}
A hyperfield $\mathbb{F}$ is called \textbf{doubly distributive} if
$$(a\boxplus b)(c\boxplus d):=\{xy|x\in a\boxplus b,y\in c\boxplus d\}=ac\boxplus ad\boxplus bc\boxplus bd$$
holds for all $a,b,c,d\in\mathbb{F}$.
\end{defn}
\begin{rem}
In general, we have $(a\boxplus b)(c\boxplus d)\subseteq ac\boxplus ad\boxplus bc\boxplus bd$.
\end{rem}
\begin{ppst}
The hyperfield of signs $\mathbb{S}$ is doubly distributive.
\end{ppst}
\begin{exam}\label{S}
Consider the polynomial $p(T)=T^3+T^2+T+1$ over the hyperfield of signs $\mathbb{S}$. One can see that
$$p(T)\in(T+1)\boxdot(T^2-T+1)=(T+1)\boxdot((T-1)\boxdot(T-1))$$
However, $0\notin\{1\}=1\boxplus 1\boxplus1\boxplus1=p(1)$ and hence $1$ is not a root of $p(T)$, which implies that $p\notin(T-1)\boxdot[(T-1)\boxdot(T+1)]$. Therefore, it is clear that
$$(T-1)\boxdot((T-1)\boxdot(T+1))\neq(T+1)\boxdot((T-1)\boxdot(T-1))$$
Hence we know that the hypermultiplication in Poly$(\mathbb{S})$ is not associative.
\end{exam}
\begin{rem}
This example shows that a hyperfield $\mathbb{F}$ is doubly distributive does not imply that the hypermultiplication in Poly$(\mathbb{F})$ is associative. Additional, as pointed out in \cite{2}, the polynomial hyperstructures over hyperfields do not always satisfy the universal property of a free algebra. In fact, there are no morphisms $\textup{Poly}(\mathbb{F})\rightarrow\mathbb{F}$ that extend the identity map $\mathbb{F}\rightarrow\mathbb{F}$ and send $T$ to the unit element in $\mathbb{F}$. One should also notice that even if $p\boxdot(q\boxdot r)\neq r\boxdot(p\boxdot q)$ for certain $p,q,r\in\textup{Poly}(\mathbb{F})$, the equality
$$p(a)\boxdot(q(a)\boxdot r(a))=r(a)\boxdot(p(a)\boxdot q(a))$$
may still true for all $a\in\mathbb{F}$, where $p(a)\boxdot q(a)=\{x\odot y:x\in p(a),y\in q(a)\}$.
\end{rem}
\begin{exam}
We have sets of polynomials over the hyperfield of signs $\mathbb{S}$
$$(T+1)\boxdot((T-1)\boxdot(T-1))=\{T^3+b_2T^2+b_1T+1:b_2,b_1\in\mathbb{S}\}$$
and
$$(T-1)\boxdot((T+1)\boxdot(T-1))=\{T^3+c_2T^2-T+1:c_2\in\mathbb{S}\}$$
In the same time, we can check that
$$(a\boxplus 1)\boxdot((a\boxplus(-1))\boxdot(a\boxplus(-1)))\boxdot(a\boxplus(-1))\boxdot((a\boxplus1)\boxdot(a\boxplus(-1)))$$
for each $a\in\mathbb{S}$.
\end{exam}

\begin{exam} Surprisingly, the hypermultiplication fails to be associative even in the polynomial hyperstructure over the simplest hyperfield: the Krasner hyperfield $\mathbb{K}$.\\
In fact, we can check that
$$(T+1)\boxdot((T^2+1)\boxdot(T+1))\neq(T^2+1)\boxdot((T+1)\boxdot(T+1))$$
in $\textup{Poly}(\mathbb{K})$.
\end{exam}

\subsection{Associativity of the hypermultiplication in Poly$(\mathbb{T})$}
Here we prove that though the hypermultiplication in Poly$(\mathbb{T})$ is associative for linear polynomials, it is not associative in general. Moreover, if $1\boxplus_{\mathbb{F}}1$ is not a singleton for $\mathbb{F}:=(\mathbb{F},\odot,\boxplus_{\mathbb{F}},1,0)$, the hypermultiplication in Poly$(\mathbb{F})$ is not associative.
\begin{lema}\label{lemma}
Given a sequence of elements $e_1\leq e_2\leq\cdots\leq e_n$ in $\mathbb{T}$, we have
$$\bigboxplus_{k=1}^ne_{i_k}=\left\{
\begin{aligned}
e_n,&&e_{n-1}< e_n\\
[-\infty,e_n],&&e_{n-1}= e_n
\end{aligned}\right.$$
where $\{i_k\}_{k=1}^n$ is a permutation of $\{1,2,\dots,n\}$.
\end{lema}

\begin{ppst}\label{T-3}
Given three elements $a,b,c$ in $\mathbb{T}$, we have
$$(0T+a)\boxdot((0T+b)\boxdot(0T+c))=\{0T^3+e_2T^2+e_1T+abc:e_2\in E_2,e_1\in E_1\}$$
where $E_2=a\boxplus b\boxplus c,E_1=ab\boxplus bc\boxplus ca$.

\begin{proof}
It is clear that
$$(0T+a)\boxdot((0T+b)\boxdot(0T+c))\subseteq\{0T^3+e_2T^2+e_1T+abc:e_2\in E_2,e_1\in E_1\}$$
Then for each $p\in(0T+a)\boxdot((0T+b)\boxdot(0T+c))$, we want to find a polynomial
$$q(T)=0T^2+d_1T+bc\in(0T+b)\boxdot(0T+c)$$
such that $p\in (0T+a)q$. In other words, we want a $d_1\in b\boxplus c$ such that
\begin{center}
$e_1\in bc\boxplus ad_1$ and $e_2\in a\boxplus d_1$
\end{center}
holds for given $e_1\in E_1$ and $e_2\in E_2$. Consider
$$d_1=\left\{\begin{aligned}
\max\{e_2,a\},&&a\leq\max\{b,c\}\\
\max\{a^{-1}e_1,a^{-1}bc\},&&a>\max\{b,c\}
\end{aligned}\right.$$
we first need to check $d_1\in b\boxplus c$. Note that when $a\leq\max\{b,c\}$, one has
$$\max\{e_2,a\}\leq\max\{a,b,c\}=\max\{b,c\}$$
so $d_1\in b\boxplus c$ in this case. When $a>\max\{b,c\}$, we have
$$\max\{a^{-1}e_1,a^{-1}bc\}\leq \max{b,c,a^{-1}bc}=\max\{b,c\}$$
so $d_1\in b\boxplus c$ in this case. Then we are going to check that
\begin{center}
$e_1\in bc\boxplus ad_1$ and $e_2\in a\boxplus d_1$
\end{center}
When $a\leq\max\{b,c\}$, one can see
$$e_1\in bc\boxplus ad_1=E_1\quad\textup{ and }\quad
e_2\in a\boxplus d_1=\left\{
\begin{aligned}
\left[-\infty,a\right],&&e_2\leq a\\
\{e_2\},&&e_2>a
\end{aligned}\right.$$
When $a>\max\{b,c\}$, one can see
$$e_1\in bc\boxplus ad_1=\left\{
\begin{aligned}
\left[-\infty,bc\right],&&e_1\leq bc\\
\{e_1\},&&e_1>bc
\end{aligned}\right.\quad\textup{ and }\quad
e_2\in a\boxplus d_1=\{a\}=E_2$$
Therefore, we are done.
\end{proof}
\end{ppst}
\begin{rem}
This proposition shows that $\boxdot$ is associative for linear polynomials in Poly$(\mathbb{T})$. In fact, we can generalize this proposition to the following one.
\end{rem}
\begin{ppst}\label{T-n}
Given a sequence of elements $\{a_i\}_{i=1}^n$ in $\mathbb{T}$, we define
$$S_k:=\bigcup_{p\in S_{k-1}}(0T+a_k)\boxdot p$$
for $k\geq2$ and $S_1=\{(0T+a_1)\}$. Then we have
$$S_n=\{0T^n+c_{n-1}T^{n-1}+\cdots+c_1T+c_0\,|\,c_{n-s}\in C_{n-s}=\bigboxplus_{I_s}\!_\mathbb{T}(\bigodot_{j=1}^s a_{t_j})\}$$
where $I_s=\{t_j\}_{j=1}^s$ represents an $s$-elements collection of $\{1,2,\dots,n\}$ and $\boxplus_{\mathbb{T}}$ represents the hyperaddition in $\mathbb{T}$.

\begin{proof} We will prove it by induction on the length of the sequence $n$.\\
It is clear that our claim is true for $n=1,2$ and Proposition \ref{T-3} tells us that our claim is also true for $n=3$.\\
Suppose our claim is true for $n\leq m$, we are going to check the case $n=m+1$.\\
With our inductive assumption, we know that
$$S_m=\{0T^m+d_{m-1}T^{m-1}+\cdots+d_1T+d_0\,|\,d_{m-s}\in D_{m-s}=\bigboxplus_{J_s}\!_\mathbb{T}(\bigodot_{j=1}^s a_{t_j})\}$$
where $J_s$ represents an $s$-elements collection of $\{1,2,\dots,m\}$. Then we are going to check that
$$\bigcup_{q\in S_m}(0T+a_{m+1})\boxdot q=\{0T^{m+1}+c_mT^{m}+\cdots+c_1T+c_0\,|\,c_{m+1-s}\in C_{m+1-s}=\bigboxplus_{I_s}\!_\mathbb{T}(\bigodot_{j=1}^s a_{t_j})\}$$
First, $D_{m-s}$ is either $\{e_s\}$ or $[-\infty,e_s]$ where $e_s=\max_{J_s}{\bigodot_{j=1}^s a_{t_j}}$ from Lemma \ref{lemma} and the second situation holds if and only if there exists two index subsets $J'_s$ and $J''_s$ of $\{1,2,\dots,m\}$ such that $\bigodot_{t'_j\in J'_s} a_{t'_j}=\bigodot_{t''_j\in J''_s} a_{t''_j}=e_s$. Similarly, $C_{m+1-s}$ is either $\{f_s\}$ or $[-\infty,f_s]$ where $f_s=\max_{I_s}\bigodot_{j=1}^s$ in $S_{m+1}$.\\
For any polynomial $p(T)=0T^{m+1}+c_{m}T^{m-1}+\cdots+c_1T+c_0$ in $(0T+a_{m+1})S_m$, we know that
$$c_{m+1-s}\in a_{m+1}d_{m-(s-1)}\boxplus_{\mathbb{T}} d_{m-s}$$
and then we are going to show that $c_i\in a_{m+1}d_{m-(s-1)}\boxplus_{\mathbb{T}}d_{m-s}\subseteq C_s$.\\
If $e_s=a_{m+1}e_{s-1}$, then it is clear that $c_{m+1-s}\in a_{m+1}d_{m-(s-1)}\boxplus_{\mathbb{T}}d_{m-s}\subseteq[-\infty,e_s]=C_s$ since we have $f_s=e_s=a_{m+1}e_{s-1}$.\\
If $e_s>a_{m+1}e_{s-1}$, then $C_s=D_s$ since $f_s=e_s$ in this case and we can not find a $J_{s-1}$ such that $a_{m+1}\odot(\bigodot_{J_{s-1}}a_{t_j})=f_s$. Therefore, $a_{m+1}d_{m-(s-1)}\boxplus_{\mathbb{T}}d_{m-s}\subseteq D_s=C_s$.\\
If $e_s<a_{m+1}e_{s-1}$, then we know that $f_s=a_{m+1}e_{s-1}$ and hence $C_s=a_{m+1}D_{s-1}$, which implies that $c_{m+1-s}\in a_{m+1}d_{m-(s-1)}\boxplus_{\mathbb{T}} d_{m-s}\subseteq a_{m+1}D_{s-1}=C_s$.\\
In conclusion, we always have $c_{m-s}\in C_s$ and hence
$$p(T)\in\{0T^{m+1}+c_mT^{m}+\cdots+c_1T+c_0\,|\,c_{m+1-s}\in C_{m+1-s}=\bigboxplus_{I_s}\!_\mathbb{T}(\bigodot_{j=1}^s a_{t_j})\} $$
which implies that
$$S_{m+1}\subseteq\{0T^{m+1}+c_mT^{m}+\cdots+c_1T+c_0\,|\,c_{m+1-s}\in C_{m+1-s}=\bigboxplus_{I_s}\!_\mathbb{T}(\bigodot_{j=1}^s a_{t_j})\}$$
Then consider a polynomial $p(T)=0T^{m+1}+c_{m}T^{m-1}+\cdots+c_1T+c_0$ with $c_{m+1-s}\in C_s$, we are going to show that there exists a $q\in S_m$ such that $p\in (0T+a_{m+1})q$.\\
First, according to our inductive assumption, we can suppose that $a_1\geq\cdots\geq a_m$ and then immediately have $e_s=\bigodot_{j=1}^sa_j$. Since we have $a_{m+1}e_{s-1}> e_s$ for $a_{m+1}> a_s$ and have $a_{m+1}e_{s-1}< e_s$ for $a_{m+1}< a_s$, there exist $s_1,s_2$ with $1\leq s_1,s_2\leq m+1$ such that
\begin{center}
$a_{m+1}e_{s-1}> e_s$ for $s> s_1$ and $a_{m+1}e_{s-1}< e_s$ for $1<s< s_2$
\end{center}
It is clear that such $s_1$ and $s_2$ depend on the order of $a_{m+1}$ in $\{a_i\}_{i=1}^{m+1}$. For example, if $a_{m+1}$ is smaller than any element of $\{a_i\}_{i=1}^m$, then $s_1=m,s_2=m+1$ and if $a_{m+1}$ is bigger than any element of $\{a_i\}_{i=1}^m$, then $s_1=1,s_2=2$.\\
We then try to find a qualifying $q(T)=0T^{m}+d_{m-1}d^{m-1}+\cdots+d_1T+d_0$. First of all, we know that $d_0=e_m$ and then want to inductively give other $d_{m-s}$.\\
If $s_1=m,s_2=1$, then we have $a_{m+1}e_{s-1}= e_s$ for any $s$. So, $a_1=a_2=\cdots=a_{m+1}$, which implies that it is a trivial case.\\
If $s_1=1,s_2=2$, then we have $a_{m+1}e_{s-1}> e_s$ for any $s$. We go from $d_0=e_m<a_{m+1}d_{m-1}$, where we have
$$d_1=\left\{
\begin{aligned}a^{-1}_{m+1}e_{m},&&c_1\in[-\infty,d_0)\\a^{-1}_{m+1}c_1,&&c_1\in[d_0,a_{m+1}e_{m-1}]\end{aligned}
\right.=\max\{a^{-1}_{m+1}d_{0},a^{-1}_{m+1}c_1\}$$
to be what we want. Generally, let $d_i=\max\{a^{-1}_{m+1}d_{i-1},a^{-1}_{m+1}c_i\}$, we have $d_i\leq e_{m-i}$. Then we can obtain $c_i\in a_{m+1}d_i\boxplus_{\mathbb{T}}d_{i-1}$ and $d_i\in D_{m-i}$ similarly. Here, we should notice that $c_m=a_{m+1}$.\\
If $s_1=m,s_2=m+1$, then we have $a_{m+1}e_{s-1}<e_s$ for any $s$. Here we do from $d_{m}=0$, where we can check
$$d_{m-1}=\left\{
\begin{aligned}c_m,&&c_m\in[a_{m+1}d_m,e_1]\\a_{m+1}d_{i+1},&&c_m\in[-\infty,a_{m+1}d_m)\end{aligned}
\right.=\max\{c_m,a_{m+1}d_{m}\}$$
is the valid one we want. In general, let $d_i=\max\{c_{i+1},a_{m+1}d_{i+1}\}$, we have $d_i\leq e_{m-i}$ and can obtain $c_i\in a_{m+1}d_i\boxplus_{\mathbb{T}}d_{i-1}$ and $d_i\in D_{m-i}$ similarly. Here $c_1=d_0=e_m>a_{m+1}d_1$.\\
In general, we can assume that $s_1>1$ and $s_2<m+1$ since those will lead to the cases solved above. Then we have exactly four cases.\\
First, when $m>s_1>1$ and $s_2=1$, we have $a_{m+1}e_{s-1}=e_s$ for $s_1\geq s>1$ and $a_{m+1}e_{s-1}>e_s$ for $s>s_1$. Then we go from $d_0=e_m$, let
\begin{center}
$d_i=\max\{a^{-1}_{m+1}d_{i-1},a_{m+1}^{-1}c_i\}$ for all $s$
\end{center}
and we can check that now the $q(T)$ is what we want just as the $s_1=1,s_2=2$ case.\\
Second, when $s_1=m$ and $1<s_2<m$, we have $a_{m+1}e_{s-1}<e_s$ for $1<s<s_2$ and $a_{m+1}e_{s-1}=e_s$ for $s\geq s_2$. Then we go from $d_m=0$, let
\begin{center}
$d_i=\max\{c_{i+1},a_{m+1}d_{i+1}\}$ for all $s$
\end{center}
and we can check that now the $q(T)$ is what we need as the $s_1=m,s_2=m+1$ case.\\
Third, when $2<s_1+1=s_2<m+1$, we have
\begin{center}
$a_{m+1}e_{s-1}<e_s$ for $1<s<s_2$ and $a_{m+1}e_{s-1}>e_s$ for $s\geq s_2$
\end{center}
Notice that in this case we have $a_{s_2}<a_{m+1}<a_{s_2-1}$, so $C_{s_2-1}=D_{s_2-1}$ and $C_{s_2}=a_{m+1}D_{s_2-1}$ are both singletons, which implies that $c_{m+1-(s_2-1)}=a_{m+1}c_{m+1-s_2}$. Now let
$$d_i=\left\{\begin{aligned}
\max\{a^{-1}_{m+1}d_{i-1},a_{m+1}^{-1}c_i\},&\quad i=0,\dots,m-s_2\\
e_{s_2-1}=c_{m+1-(s_2-1)},&\quad i=m-(s_2-1)\\
\max\{c_{i+1},a_{m+1}d_{i+1}\},&\quad i=m-(s_2-2),\dots,m-1
\end{aligned}\right.$$
Then we can clearly see $c_i\in a_{m+1}d_i\boxplus_{\mathbb{T}}d_{i-1}$ and $d_i\in D_{m-i}$ in this case.\\
Forth, when $2<s_1+1<s_2<m+1$, we have
\begin{center}
$a_{m+1}e_{s-1}<e_s$ for $1<s<s_2$,\\
$a_{m+1}e_{s-1}=e_s$ for $s_1\geq s\geq s_2$,\quad $s>s_1$ for $a_{m+1}e_{s-1}>e_s$.
\end{center}
Note that here we have $a_{s_1+1}<a_{m+1}<a_{s_2-1}$, so $C_{s_2-1}=D_{s_2-1}$ and $C_{s_1+1}=a_{m+1}D_{s_1}$ are both singletons, which implies that $c_{m+1-(s_2-1)}=e_{s_2-1}$ and $c_{m+1-(s_1+1)}=a_{m+1}e_{s_1}$. Now let
$$d_i=\left\{\begin{aligned}
\max\{a^{-1}_{m+1}d_{i-1},a_{m+1}^{-1}c_i\},&\quad i=0,\dots,m-(s_1+1)\\
e_{m-i},&\quad i=m-s_1,\dots,m-(s_2-1)\\
\max\{c_{i+1},a_{m+1}d_{i+1}\},&\quad i=m-(s_2-2),\dots,m-1
\end{aligned}\right.$$
Then one can clearly see $c_i\in a_{m+1}d_i\boxplus_{\mathbb{T}}d_{i-1}$ and $d_i\in D_{m-i}$ in this case.\\
In conclusion, we can always find such $q$ with $p\in(T+a_{m+1})q$, which follows that
$$S_{m+1}\supseteq\{0T^{m+1}+c_mT^{m}+\cdots+c_1T+c_0\,|\,c_{m+1-s}\in C_{m+1-s}=\bigboxplus_{I_s}\!_\mathbb{T}(\bigodot_{j=1}^s a_{t_j})\}$$
and then our claim is true for $n=m+1$. Therefore, we are done.
\end{proof}
\end{ppst}
\begin{rem}
The proof shows that we can always find a 'solution' $q\in S_m$ such that $p\in(T+a_{m+1})q$ for each $p\in S_{m+1}$. One could notice that the examples in the non-associativity part imply that we cannot find such a 'solution' even for some really simple cases. Following is the theorem 4.1 in \cite{2} which is very important in polynomial algebra over the tropical hyperfield $\mathbb{T}$.
\end{rem}
\begin{thm} Given a monic polynomial $p(T)$ of degree $n$ in \textup{Poly}$(\mathbb{T})$, then\\
(1) There is a unique sequence $a_1,\dots,a_n\in\mathbb{T}$, up to permutation of indices, such that
$$p\in \{0T^n+c_{n-1}T^{n-1}+\cdots+c_1T+c_0\,|\,c_{n-s}\in C_{n-s}=\bigboxplus_{I_s}\!_\mathbb{T}(\bigodot_{j=1}^s a_{t_j})\}$$
where $I_s=\{t_j\}_{j=1}^s$ represents an $s$-elements collection of $\{1,2,\dots,n\}$ and $\boxplus_{\mathbb{T}}$ represents the hyperaddition in $\mathbb{T}$.\\
(2) The equalities $\textup{mult}_a(p)=\textup{Card}(\{i\in\{1,\dots,n\}\,|\,a=a_i\})$ hold for all $a\in\mathbb{T}$.
\end{thm}
\begin{rem}
This theorem and Proposition \ref{T-n} show that each monic polynomial over $\mathbb{T}$ belongs to a unique hyperproduct (up to permutation) of linear polynomials and hyperproducts of polynomials in Poly$(\mathbb{T})$ that can be exactly linear-represented are also unique up to permutation. However, as Oliver Lorscheid points out, for the associativity, we still need to consider the polynomials in Poly$(\mathbb{T})$ which cannot be represented as a hyperproduct.
\end{rem}
\begin{defn}
A polynomial $p$ over a hyperfield $\mathbb{F}$ is \textbf{reducible}, if it is exactly the hyperproduct of  two positive-degreed polynomials. Otherwise, it is \textbf{irreducible}.
\end{defn}
\begin{exam}
The polynomial $p(T)=0T^2+2$ over the tropical hyperfield $\mathbb{T}$ is irreducible. Otherwise, if it is reducible, then there exist two linear polynomials $q_1(T)=a_1T+b_1$ and $q_2(T)=a_2T+b_2$ with $q_1(T)\boxdot q_2(T)=\{0T^2+2\}$. In this case, one has $a_1b_2\boxplus a_2b_1=\{-\infty\}$. Since $a_1,a_2$ can not be the zero element $-\infty$, we know $b_1=b_2=-\infty$ and then have $-\infty=b_1b_2=2$, a contradiction.
\end{exam}

\begin{thm}
If the hypermultiplication is associative for all reducible polynomials over the tropical hyperfield $\mathbb{T}$, then this hypermultiplication is associative for all polynomials in $\textup{Poly}(\mathbb{T})$.
\end{thm}

However, the tropical hyperfield cannot meet the condition. In fact, we have the following observation.
\begin{thm}
Given a hyperfield $(\mathbb{F},\odot,\boxplus_{\mathbb{F}},1,0)$, if $1\boxplus1$ is not a singleton, then the hypermultiplication $\mathbb{F}$ in $\textup{Poly}(\mathbb{F})$ is not associative.

\begin{proof}
We can see the following sets
$$\begin{aligned}
&(T^2+1)\boxdot((T+1)\boxdot(T+1))\\
&=\{(T^2+1)\boxdot(T^2+dT+1)\,|\,d\in1\boxplus_{\mathbb{F}}1\}\\
&=\{T^4+dT^3+(1\boxplus_{\mathbb{F}}1)T^2+dT+1\,|\,d\in1\boxplus_{\mathbb{F}}1\}
\end{aligned}$$
and
$$\begin{aligned}
&(T+1)\boxdot((T^2+1)\boxdot(T+1))\\
&=(T+1)\boxdot(T^3+T^2+T+1)\\
&=\{T^4+d_1T^3+(1\boxplus_{\mathbb{F}}1)T^2+d_2T+1\,|\,d_1,d_2\in1\boxplus_{\mathbb{F}}1\}
\end{aligned}$$
are not equal.
\end{proof}
\end{thm}


\end{document}